\documentclass[a4paper,10pt,final]{amsart}
\usepackage[english]{certus}
\usepackage{amsmath}
\usepackage{flexisym,breqn}
\usepackage[notref,notcite,draft]{showkeys} 
\usepackage{verbatim} 

\newcommand{\NN}{{\mathbb N}}
\newcommand{\ZZ}{{\mathbb Z}}

\newcommand{\CC}{{\mathbb C}}

\newcommand{\MM}{{\mathbb M}}
\newcommand{\GG}{{\mathbb G}}

\newcommand{\A}{{\mathcal A}}
\DeclareMathOperator{\card}{card}

\renewcommand{\max}{{\operatorname{max}}}
\newcommand{\varps}{{\varepsilon}}
\newcommand{\rg}{{\operatorname{rg\hspace{0.04cm}}}}
\newcommand{\tens}{\otimes}
\newcommand{\spann}{{\operatorname{span}}}
\newcommand{\To}{\longrightarrow}

\newcommand{\Tor}{\operatorname{Tor}}
\newcommand{\op}{{\operatorname{{op}}}}
\newcommand{\bet}{\beta^{(2)}}
\newcommand{\Pol}{{\operatorname{Pol}}}
\newcommand{{\hs}}{{\operatorname{HS}}}
\newcommand{\alg}{{\operatorname{alg}}}

\newcommand{\st}{{\operatorname{st}}}

\newcommand{\vN}{{\operatorname{vN}}}

\author{Vadim Alekseev} 
\address{Vadim Alekseev,
Mathematisches Institut,
Georg-Au\-gust-Uni\-versi\-t{\"a}t G{\"o}t\-ting\-en,
Bunsenstra{\ss}e 3-5,
D-37073 G{\"o}ttingen, 
Germany.}
\email{alekseev@uni-math.gwdg.de}

\author{David Kyed} 
\address{David Kyed,
Department of Mathematics,
KU Leuven,
Celestijnenlaan 200B,
B-3001 Leuven, 
Belgium.}
\email{David.Kyed@wis.kuleuven.be}

\title[Amenability and  vanishing of $L^2$-Betti numbers]{Amenability and  vanishing of $L^2$-Betti numbers: \\ an operator algebraic approach}

\keywords{Amenability, $L^2$-Betti numbers, operator algebras}
\subjclass[2010]{46L10, 43A07, 18G15}

\thanks{The research of the second named author is supported by The Danish Council for Independent Research $|$ Natural Sciences and  the ERC Starting Grant VNALG-200749}

\begin{document}

\begin{abstract}
We introduce a F{\o}lner condition for dense subalgebras in finite von Neumann algebras and prove that it implies dimension flatness of the inclusion in question. It is furthermore proved that the F{\o}lner condition naturally generalizes  the existing notions of amenability and that the ambient von Neumann algebra of a F{\o}lner algebra is automatically injective.
As an application, we show how our techniques unify previously known results  concerning vanishing of $L^2$-Betti numbers for amenable groups, quantum groups and groupoids and moreover provide a large class of new examples of algebras with vanishing $L^2$-Betti numbers.
\end{abstract}

\maketitle

\section{Introduction} 
$L^2$-Betti numbers originated from topology, but have later unveiled an interesting connection to operator algebras. Indeed, while the original definition of  $L^2$-Betti numbers for groups involved a geometric construction \cite{atiyah-l2,cheeger-gromov}, it was later shown by L\"uck (see also \cite{farber-vN-categories}) that they can also be described in terms of certain $\Tor$ modules  ---  thus using only homological algebra and modules over von Neumann algebras. This has led to several generalizations of $L^2$-Betti numbers: they have been defined for discrete measured groupoids \cite{sauer-betti-for-groupoids}, quantum groups \cite{quantum-betti} and general subalgebras of finite von Neumann algebras \cite{CS}. Similarly, $L^2$-Betti numbers were introduced in the setting of equivalence relations in \cite{gaboriau-L2-for-eq-rel} following the original geometric approach, but as shown in \cite{sauer-betti-for-groupoids} these can also be expressed in terms of homological algebra.

There are two properties which attract attention in all these situations. Firstly, in each case the definition involves certain operator algebras canonically associated with the situation at hand. Secondly, for groups, quantum groups and groupoids it is well known that the $L^2$-Betti numbers vanish in the  presence of amenability. In view of this, it is natural to seek a common operator algebraic reason for this to happen. In doing so, one firstly observes that the actual reason for the vanishing of $L^2$-Betti numbers for amenable groups and quantum groups is a certain dimension-flatness property of an inclusion $\mc A\subseteq M$, where $\mc A$ is a strongly dense $\ast$-subalgebra of a finite von Neumann algebra $M$. Secondly, the key to the proof of this dimension-flatness result is a F{\o}lner condition for the notion of amenability in question.
In the present paper we introduce a F{\o}lner type condition for a general weakly dense $*$-subalgebra $\mc A$ in a tracial von Neumann algebra $(M,\tau)$ and show  how this leads to dimension-flatness  of the inclusion $\mc A\subseteq M$ and subsequently to the vanishing of the operator algebraic $L^2$-Betti numbers $\beta_p^{(2)}(\mc A, \tau)$. 
This approach unifies the above mentioned vanishing results and furthermore provides a large class of new examples of algebras with vanishing $L^2$-Betti numbers (see Section \ref{example-section}). More precisely, we prove the following:

\begin{thm*}[see Theorem \ref{dim-flat-thm} \& Corollary \ref{CS-vanishing}]
If $\mc A \subseteq M$ satisfies the F{\o}lner condition then for any left $\mc A$-module $X$ and any $p\geqslant 1$ we have
\[
\dim_M \Tor_p^{\mc A}(M,X)=0,
\]
and the Connes-Shlyakhtenko $L^2$-Betti numbers of $\mc A$ vanish in positive degrees.
\end{thm*}

Secondly, we link our F{\o}lner condition to the classical operator algebraic notion of amenability by proving (a slightly more general version of) the following:

\begin{thm*}[see Theorem \ref{injectivity-thm}] If $\mc A\subseteq M$ satisfies the F{\o}lner condition $M$ is injective.
\end{thm*} 

The  F\o{}lner condition consists of two requirements: an almost invariance property and a trace approximation property. 
The  almost invariance property requires that the left action of $\mc A$ on $L^2(M)$ admits almost invariant subspaces, the almost invariance being measured by a dimension function. The naive approach would be to use the usual dimension over $\mb C$ and require these subspaces to be finite-dimensional; however, it turns out that the dimension theory over a von Neumann subalgebra $N\subseteq \mc A$ can be efficiently applied in order to allow for almost invariant subspaces which are finitely generated projective  $N$-modules (thus usually having infinite dimension over $\mb C$). For this reason, parts of the paper deals with results related to dimension theory for modules over von Neumann algebras; the reader not familiar with these notions may, without loosing the essential ideas, think of the case $N=\CC$, which reduces most arguments to finite dimensional linear algebra. \\

\noindent{\emph{Structure:}} 
The paper is organized as follows. In the second section we recapitulate L\"uck's dimension theory for modules over von Neumann algebras and prove a few results concerning a relative dimension function in this context. In the third section we introduce the operator-algebraic F\o{}lner condition and draw some easy consequences, preparing for the proof of the main theorem which is given in the fourth section. The fifth section is devoted to the discussion concerning the relationship between the F\o{}lner condition for an algebra and injectivity of  its enveloping  von Neumann algebra. In the sixth section we discuss examples and show how the F{\o}lner condition implies vanishing of $L^2$-Betti numbers in a variety of different instances. \\

\noindent{\emph{Assumptions:}} Throughout the paper, all generic von Neumann algebras are assumed to have separable predual. Our main focus will be on finite von Neumann algebras, for which separability of the predual is equivalent to separability of the von Neumann algebra itself for the strong, weak and ultra-weak topology. These topologies are furthermore all metrizable on bounded subsets of $M$.  In order to be consistent with the  general Hilbert $C^*$-module terminology (see below) inner products on Hilbert spaces are assumed to be linear in the second variable and conjugate linear in the first. Algebraic tensor products are denoted ``$\odot$'', tensor products of Hilbert spaces by ``$\tens$'' and tensor products of von Neumann algebras by $``\bar{\tens}$''. Our main setup  will consist of an inclusion $N\subseteq M$ of finite von Neumann algebras with a fixed normal tracial state $\tau$ on $M$. The von Neumann algebra $N$ should be thought of as a ``coefficient algebra'' and will in many applications be equal to $\CC$.  We denote by $H$ the GNS-space $L^2(M,\tau)$, by $J\colon H\to H$ the modular conjugation arising from $\tau$ and by $E\colon M\to N$ the $\tau$-preserving conditional expectation. Note that $H$ automatically carries a normal right $N$-action given by $\xi\cdot x:= Jx^*J\xi$ which extends the natural right $N$-module structure on $M$. \\

\noindent{\emph{Acknowledgements:}}  The authors would like to thank Thomas Schick for a crucial remark at an early stage in the project and Stefaan Vaes for suggesting an improvement of Theorem \ref{injectivity-thm}. Moreover, thanks are due to {\'{E}}tienne Blanchard, Henrik D.~Petersen and the anonymous referee for numerous valuable corrections and suggestions.

\section{Dimension theory}\label{dim-theory-section}
In this section we recapitulate parts of the von Neumann dimension theory  for modules over a finite von Neumann algebra  and introduce a relative version of this notion of dimension. The relative dimension function will play a prominent role in the definition of the F{\o}lner condition given in Section \ref{foelner-section}.\\ 

\subsection{Dimension of modules and the trace on endomorphisms}
Consider a finite von Neumann algebra $N$ endowed with a normal, faithful tracial state $\tau$ and denote by $L^2(N)$ the GNS-space arising from $\tau$ and by $N^{\op}$ the opposite von Neumann algebra.  Denoting by $\Lambda$ the natural embedding of $N$ into $L^2(N)$ the inner product in $L^2(N)$ is therefore given by $\ip{\Lambda(x),\Lambda(y)}=\tau(x^*y)$ for $x,y\in N$.  In what follows, we will often suppress the map $\Lambda$ and simply consider $N$ as a subspace of $L^2(N)$. A Hilbert space $F$ endowed with a normal $*$-representation of $N^{\op}$ is called a \emph{finitely generated normal (right) $N$-module} if there exists an isometric $N^{\op}$-equivariant embedding of $F$ into $L^2(N)^k$ for some finite integer $k$, where $L^2(N)^k$ is considered with the standard $N^{\op}$-action given by
\[
a^{\op}\colon(\xi_1,\dots, \xi_n)\mapsto (Ja^*J\xi_1,\dots, Ja^*J\xi_n).
\]
Thus, any finitely generated normal $N$-module $F$ is isomorphic to one of the form $pL^2(N)^k$ for a projection $p\in \MM_k(N)$ and \emph{the von Neumann dimension} of $F$ is then defined as
\[
\dim_N F:=\sum_{i=1}^k\tau(p_{ii}).
\]
The number $\dim_N F$ is  independent of the choice of $p$ and $k$ (see e.g.~\cite[Section 1.1.2]{Luck02}) and thus depends only on the isomorphism class of $F$. 
\begin{Rem}
In most of the existing literature,  finitely generated normal right modules over $N$ are called \emph{finitely generated Hilbert $N$-modules}, but since we are going to consider finitely generated Hilbert $C^*$-modules over $N$ (in the operator algebraic sense \cite{lance}) as well, we have chosen the term ``normal module'' in order to avoid unnecessary confusion. We advise the reader not familiar with these notions to consult the article \cite{frank-hilbertian-vs-hilbert} where a detailed comparison is carried out.
\end{Rem}

Next we recall the construction of the trace on the endomorphism von Neumann algebra $\End_N(F)$ associated with a finitely generated normal $N$-module $F$; here, and in what follows, $\End_N(F)$ denotes the von Neumann algebra of all bounded operators on $F$ commuting with the $N$-action. The trace on $\End_N(F)$  was previously considered by L{\"u}ck \cite[Section 1.1]{Luck02}  and also by Farber \cite{farber-vN-categories} in the more general context of von Neumann categories. 

\begin{Lemma}[{\cite[Section 1.1.3]{Luck02}}]\label{trace-lemma}
Let $F$ be a finitely generated normal $N$-mo\-dule and consider an operator $T\in \End_N(F)$. Upon choosing an isomorphism $F\simeq pL^2(N)^k$ we obtain an induced $*$-isomorphism  $\End_N(F)\simeq p\MM_k(N)p$ and we denote the matrix representing $T$ by $(T_{ij})_{i,j=1}^k$.  Then the number $\Tr_{N}(T):= \sum_{i=1}^{k}\tau(T_{ii})$ does not depend on the choice of isomorphism and defines a normal, faithful, positive trace on $\End_N(F)$ with $\Tr_N(\id_F)=\dim_NF$. 
\end{Lemma}

\begin{Rem}
In the case $N=\CC$ the normal $N$-module $F$ is simply a finite dimensional vector space and hence $\Tr_N(-)$ is just the standard, non-normalized trace on $B(F)$.  
\end{Rem}

The first step towards extending the notion of dimension to arbitrary $N$-modules is to pass from normal modules to Hilbert $C^*$-modules over $N$.

\begin{Def}
A (right) Hilbert $C^*$-module over $N$ consists of an algebraic right $N$-module $\mc X$ together with a map $\ip{\cdot,\cdot}_N\colon \mc X\times \mc X\to N$ such that
\begin{itemize}
\item[(i)] For all $x\in \mc X$ we have ${\ip{x,x}}_N\geqslant 0$ and ${\ip{x,x}}_N=0$ only for $x=0$.
\item[(ii)] For all $x,y\in \mc X$ we have $\ip{x,y}_N=\ip{y,x}_N^*$.
\item[(iii)] For all $x,y,z\in \mc X$ and $a,b\in N$ we have $\ip{z,xa+yb}_N= \ip{z,x}_Na +\ip{z,y}_Nb$.
\item[(iv)] The space $\mc X$ is complete with respect to the norm $\|x\|:=\|\ip{x,x}\|^{\frac12}$.
\end{itemize}
A Hilbert $C^*$-module over $N$ is called \emph{finitely generated} if it is algebraically finitely generated over $N$ and \emph{projective} if it is projective as a right $N$-module.  Moreover, the terms ``Hilbert module over $N$'', ``Hilbert $C^*$-module over $N$'' and ``Hilbert $N$-module'' will be used synonymously in the sequel.
\end{Def}

Note that the finitely generated free $N$-modules $N^k$ become Hilbert $N$-modules when endowed with the $N$-valued inner product $\ip{x,y}_{\st}:=\sum_{i=1}^k x_i^*y_i$; we refer to this as the \emph{standard} inner product on $N^k$.  If $(\mc P, \ip{\cdot,\cdot}_N)$ is a finitely generated projective Hilbert $N$-module then $\tau\circ\ip{\cdot,\cdot}_N$ defines a positive definite $\CC$-valued inner product on $\mc P$ and we denote the Hilbert space completion of $\mc P$ with respect to this inner product by $L^2(\mc P)$. It was shown by L\"uck that we get a finitely generated normal $N$-module in this way and that this construction yields an equivalence of  categories:

\begin{Thm}[{L{\"u}ck, \cite[Lemma 6.23 \& Theorem 6.24]{Luck02}}]\label{lueck-eq-of-categories}
Every finitely generated projective $N$-module allows an $N$-valued inner product turning it into a Hilbert $N$-module and the inner product is unique up to unitary isomorphism. Furthermore, Hilbert space completion constitutes an equivalence of categories between the category of finitely generated projective Hilbert $N$-modules and the category of finitely generated normal $N$-modules. 
\end{Thm}
We remark that a finitely generated projective Hilbert $N$-module $(\mc P, \ip{\cdot, \cdot}_N)$ is automatically self-dual; i.e.~$\Hom_N(\mc P,N)=\{\ip{x, \cdot}_N \mid x\in \mc P\}$. This follows easily from the uniqueness of the inner product and the obvious self-duality of the finitely generated free Hilbert modules $(N^k,\ip{\cdot,\cdot}_{\st})$. Furthermore, by \cite[Lemma 2.3.7]{manuilov-troitsky} every finitely generated Hilbert submodule in $\mc P$ splits off orthogonally as a direct summand and is therefore, in particular, itself finitely generated and projective.  Due to Theorem \ref{lueck-eq-of-categories}, it makes sense, for a finitely generated projective $N$-module $\mc P$,  to 
define the von Neumann dimension of  $\mc P$ as $\dim_N \mc P=\dim_N L^2(\mc P)$ where $L^2(\mc P)$ is the Hilbert space completion relative to some choice of Hilbert $N$-module structure on $\mc P$. Moreover, for an arbitrary (algebraic) $N$-module $\mc X$ its von Neumann dimension is defined as
\[
\dim_N\mc X = \sup\{\dim_N\mc P \mid \mc P\subseteq \mc X, \mc P \ \text{finitely generated projective}  \}\in [0,\infty].
\]
This dimension function is no longer faithful (i.e.~non-zero modules might be zero-dimensional) but it still has remarkably good properties; one of its prime features being  additivity with respect to short exact sequences \cite[Theorem 6.7]{Luck02}.  For a submodule $\mc E$ in an $N$-module $\mc F$ the \emph{algebraic closure} $ \overline{\mc E}^{\alg} $ of $\mc E$ inside $\mc F$ is defined as the intersection of the kernels of those homomorphisms in the dual module $\Hom_N(\mc F,N)$ that vanish on $\mc E$, and if $\mc F$ is finitely generated we have  that
$
\dim_N \mc E = \dim_N  \overline{\mc E}^{\alg}
$ 
\cite[Theorem 6.7]{Luck02}.

\begin{Lemma}\label{same-closure-lem}
If $\mc P$ is a finitely generated projective Hilbert $N$-module and $\mc E\subseteq \mc P$ a submodule then ${\mc E}^{\perp\perp}=\overline{\mc E}^{\alg}$  
and $\mc E^{\perp\perp}$ splits off orthogonally as a direct summand in $\mc P$. Furthermore, the Hilbert space closures of $\mc E$ and $\mc E^{\perp\perp}$ in $L^2(\mc P)$ coincide.
\end{Lemma}
\begin{proof}
Since $\mc P$ is selfdual, the equality $\mc{E}^{\perp\perp}=\overline{\mc E}^\alg$ follows directly from the definition of the algebraic closure. Since $\mc P$ is finitely generated and projective it now follows from \cite[Theorem 6.7]{Luck02}  that the same is true for $\mc E^{\perp\perp}$ and by \cite[Lemma 2.3.7]{manuilov-troitsky} it therefore follows that $\mc{E}^{\perp\perp}$ splits off as an orthogonal direct summand. Since $\dim_N (\mc E^{\perp\perp}/\mc E) =0$,  Sauer's local criterion \cite[Theorem 2.4]{sauer-betti-for-groupoids} implies that for every $x\in \mc E^{\perp\perp}$ there exists a sequence of projections $p_i\in N$ such that $\lim_i\tau(p_i)=1$ and $xp_i \in \mc E$. Thus
\begin{align*}
\|x-xp_i \|_2^2 &= \tau\left(\ip{x(1-p_i),x(1-p_i)}_N\right)\\
&=\tau\left(\ip{x,x}_N(1-p_i)^2\right)\\
&\leqslant \tau\left(\ip{x,x}_N^2\right)^{\frac12} \tau\left(1-p_i\right)^{\frac12}\underset{i\to\infty}{\To}0.
\end{align*}
Hence $\mc E^{\perp\perp} \subseteq \overline{\mc E}^{\|\cdot\|_2}$ and since $\mc E \subseteq \mc E^{\perp\perp}$ this proves that the two closures coincide.

\end{proof}

If $(\mc F, \ip{\cdot, \cdot}_N)$ is a finitely generated projective Hilbert $N$-module we will often identify the algebraic endomorphism ring $\End_N(\mc F)$ with $\End_N(L^2(\mc F))$ and thus also consider $\Tr_{N}(-)$ as defined on $\End_N(\mc F)$. Recall,  e.g.~from \cite[Lemma 1.7]{kajiwara-watatani}, that every finitely generated projective Hilbert $N$-module $(\mc F,\ip{\cdot,\cdot}_N)$ has a \emph{basis}; that is there exist $u_1,\dots, u_n\in \mc F$ such that $x=\sum_{i=1}^n u_i\ip{u_i,x}_N$ for each $x\in \mc F$.\footnote{Note that in general the elements $u_1,\dots,u_n$ are not linearly independent over $N$ and thus not a basis in the standard sense of homological algebra.} Furthermore, the matrix $p=(\ip{u_i,u_j}_N)_{i,j=1}^n$ is a projection in $\MM_n(N)$ and the map $\alpha\colon (\mc F,\ip{\cdot, \cdot}_N) \to (pN^{n},\ip{\cdot,\cdot}_\st)$ given by $\alpha(x)=(\ip{u_i,x}_N)_{i=1}^n$ is a unitary isomorphism of Hilbert $N$-modules with $u_i=\alpha^{-1}(pe_i)$, where $e_1,\dots,e_n$ is the standard basis in $N^n$. From this we obtain the following result.

\begin{Lemma}\label{basis-lemma}
For a finitely generated projective Hilbert $N$-module $(\mc F,\ip{\cdot,\cdot}_N)$ and an endomorphism $T\in \End_{N}(\mc F)$   we have  $ \Tr_N(T)= \sum_{i=1}^n \tau( \ip{u_i,Tu_i}_N)$ for any basis $u_1,\dots, u_n$ in $\mc F$.
\end{Lemma}

\subsection{The relative dimension function}
Next we consider a trace-preserving inclusion of $N$ into a bigger finite von Neumann algebra $(M,\tau)$ acting on its GNS-space $H=L^2(M,\tau)$. Note that any $N$-submodule $\mc F\subseteq M$ acquires an $N$-valued inner product arising from the trace-preserving conditional expectation $E\colon M\to N$ by setting $\ip{a,b}_N:=E(a^*b)$. Consider now a finitely generated, projective $N$-submodule $\mc F \subseteq M$   which is complete with respect to this pre-Hilbert module structure.  
In other words, the relation $\ip{a, b}_N:=E(a^*b)$ defines a Hilbert $N$-module structure on $\mc F$  and hence $L^2(\mc F)$ can be obtained by completion with respect to the $\CC$-valued inner product $\tau\circ \ip{\cdot,\cdot}_N$.  Since $E$ is trace-preserving, this completion is exactly  the closure of $\mc F$ in $H$. 

\begin{Def}\label{complete-defi}
An $N$-submodule $\mc F \subseteq M$ which is a Hilbert $N$-module with respect to the $N$-valued inner product arising from the conditional expectation $E\colon M\to N$ is called complete. 
\end{Def}

Let $\mc F\subseteq M$ be a complete, non-zero, finitely generated projective $N$-submodule  and denote by $P_F\in B(H)$ the projection onto its closure $F$ in $H$. Note that $F$ is a finitely generated normal $N$-module. Any operator $T\in (JNJ)'\subseteq B(H)$ gives rise to an $N$-equivariant operator $P_F TP_F|_{F}$ on $F$ and we may therefore define a normal state $\varphi_{\mc F}$ on $(JNJ)'$ by setting
\begin{align}\label{phi-F-def}
\varphi_{\mc F}(T)=\frac{\Tr_{N}(P_FTP_F|_{F} )}{\dim_N(\mc F)}.
\end{align}
By choosing a basis $u_1,\dots, u_k\in \mc F$ for the Hilbert $N$-module structure arising from $E$, the state $\varphi_{\mc F}$ can be computed as $\varphi_{\mc F}(T)=(\dim_N\mc F)^{-1}\sum_{i=1}^k \ip{\Lambda(u_i),T\Lambda(u_i)}$ where the inner product is taken in $H$. This is due to the fact that $P_FTP_F|_{\mc F}\in \End_N(\mc F)$ and hence by Lemma \ref{basis-lemma}   
\begin{align*}
\varphi_{\mc F}(T)&=\frac{\sum_{i=1}^k \tau( \ip{u_i  , P_FTP_{F} u_i }_N )}{\dim_N \mc F} \\
&= \frac{\sum_{i=1}^k \tau( E(  u_i^*(P_FTP_F u_i))}{\dim_N \mc F}\\
&=\frac{\sum_{i=1}^k \ip{\Lambda (u_i), P_FTP_F \Lambda (u_i) }}{\dim_N \mc F}\\
&=\frac{\sum_{i=1}^k \ip{\Lambda(u_i), T \Lambda( u_i) }}{\dim_N \mc F}.
\end{align*}

In what follows we will also consider operators acting on a $k$-fold amplification $H^k=\bigoplus_1^k H$ of the Hilbert space $H$; these amplifications are always implicitly assumed equipped with the diagonal normal right action of $N$. For a subspace $L\subseteq H$ we denote by $L^k\subseteq H^k$ its amplification.

\begin{Def}\label{rel-dim-defi}
Let $k\in \NN$ and  a complete, non-zero, finitely generated, projective right $N$-submodule $\mc F$ in $M$ be given.  For an $N$-submodule $\mc E\subseteq H^k$ we define \emph{the dimension of $\mc E$ relative to $\mc F$} as 
\[
\dim_{\mc F}(\mc E)=\frac{\Tr_{N}(P_{F^k}P_EP_{F^k}|_{F^k})}{\dim_N\mc F}, 
\]
where $F$ and $E$ are the  closures in $H^k$ of $\mc F$ and $\mc E$, respectively.
\end{Def}
Note that this is well-defined since $P_E$ commutes with the right action of $N$ so that $P_{F^k}P_{E}P_{F^k}|_{ F^k}\in \End_N(F^k)$. Note also the trivial fact that if two submodules have the same closure in $H^k$ then their relative dimensions agree.

\begin{Prop}\label{rel-dim-properties-amplified}
The relative dimension function $\dim_{\mc F}(-)$ has the following pro\-perties.
\begin{itemize}
\item[(i)] If $ \mc E_1\subseteq  \mc E_2\subseteq H^k$ are two  $N$-submodules then $\dim_{ \mc F}(\mc  E_1)\leqslant \dim_{ \mc F}(\mc  E_2)$.
\item[(ii)] If $ \mc E\subseteq \mc  F^k$  then $\dim_{\mc  F}(\mc  E)= \frac{\dim_N \mc  E}{\dim_N \mc F}$.  In other words, $\dim_N\mc E=\Tr_N(P_E|_{F^k})$ where $P_E$ is the projection onto the closure $E$ of $\mc E$ in $H^k$.
\end{itemize}   
\end{Prop}

\begin{proof}   
For the sake of notational convenience we restrict attention to the case $k=1$, but the same proof goes through in higher dimensions.
The first claim follows directly from positivity of the state $\varphi_{\mc F}$.  To prove the second claim, first note that $\mc F$ splits orthogonally as $\mc E^{\perp\perp}\oplus \mc G$ for some Hilbert $N$-submodule $\mc G$.  Since $\mc F$ is assumed to be complete, its $L^2$-completion coincides with its closure $F$ inside $H$ and hence 
\[
F=\overline{\mc E ^{\perp\perp}}\oplus \overline{\mc G}
\]
as an orthogonal direct sum inside $H$ and $L^2(\mc E^{\perp\perp})=\overline{\mc E ^{\perp\perp}}$. Denote by $P_{\mc E^{\perp\perp}}\in \End_N(\mc F)$ the projection onto the summand ${\mc E^{\perp\perp}}$ and by $P_{\mc E^{\perp\perp}}^{(2)}\in \End_N(F)$ its extension to $F$. Then clearly $P_{\mc E^{\perp\perp}}^{(2)}$ projects onto the summand $\overline{\mc E^{\perp\perp}}$ which by Lemma \ref{same-closure-lem} coincides with $\overline{\mc E}$. Thus $P_FP_EP_F|_{F}=P_{\mc E^{\perp\perp}}^{(2)}$. By choosing a basis  $u_1,\dots u_n$ for $\mc{E}^{\perp\perp}$ and a basis $v_1,\dots, v_l$ for $\mc G$ we get a basis for $\mc F$ and using this basis to compute the trace we get
\begin{align*}
\dim_{\mc F}\mc E &=\frac{\sum_{i=1}^n \ip{u_i, P_{\mc E^{\perp\perp}}^{(2)} u_i} +\sum_{j=1}^l \ip{v_j, P_{\mc E^{\perp\perp}}^{(2)}v_j} }{\dim_N \mc F}\\
&=\frac{\sum_{i=1}^n\ip{u_i, u_i}}{\dim_N\mc F}\\
&=\frac{\dim_N\mc E^{\perp\perp}}{\dim_N \mc F}\\
&=\frac{\dim_N \mc E}{\dim_N \mc F},
\end{align*}
where the last equality follows from Lemma \ref{same-closure-lem}.

\end{proof}
Relative dimension functions were originally introduced by Eckmann \cite{eckmann} in a topological setting (see also \cite{elek-zdc}) and the relative dimension function $\dim_{\mc F}(-)$ introduced above can be seen as an operator algebraic analogue of this construction. A similar construction for quantum groups was considered in \cite{kyed-thom}. We end this section with a small lemma which will turn out useful in the following.
\begin{Lemma}\label{intersection-lem}
 Let $ \mc  F$ be a non-zero, finitely generated, projective $N$-module and let $\mc  S_1, \dots,\mc  S_n$ be a family of submodules in  $\mc F $. Suppose that
\beqn
\forall i\in\{1,\dots,n\}:\quad \frac{\dim_N \mc  S_i }{\dim_N \mc  F} > 1 - \eps_i
\eeqn
for some $\eps_1,\dots, \eps_n>0$. Then
\beqn
\frac{\dim_N \left(\bigcap_{i=1}^n\mc  S_i\right) }{\dim_N \mc F} > 1 - \sum_{i=1}^n\eps_i.
\eeqn
\end{Lemma}
\begin{proof}
By induction it is sufficient to consider the case $n=2$. In this case we have an exact sequence of $N$-modules
\[
0 \To \mc S_1\cap \mc S_2 \To \mc S_1\oplus \mc S_2 \To {\mc S_1+\mc S_2}\To 0,
\]
and by \cite[Theorem 6.7]{Luck02}  we have $\dim_N(\mc S_1\oplus \mc S_2)=\dim_N(\mc S_1\cap \mc S_2)+\dim_N({\mc S_1+\mc S_2})$. Thus
 
\beqn
\frac{\dim_N( \mc S_1 \cap \mc S_2)}{\dim_N \mc  F} = \frac{\dim_N \mc  S_1 + \dim_N \mc S_2 - \dim_N({ \mc S_1+ \mc S_2})}{\dim_N  \mc F} > 1 - \eps_1 - \eps_2.
\eeqn
\end{proof}
Of course the whole theory has a mirror counterpart for \emph{left} modules over $N$ and we will use both theories without further specification throughout the paper. Moreover, we will not make any notational distinction between the corresponding dimension functions; so when $\mc X$ is a \emph{left} $N$-module $\dim_N\mc X$ will also denote the von Neumann dimension of $\mc X$ as a left module and similarly with the relative dimension functions.

\section{The F{\o}lner property}\label{foelner-section}
Consider again a trace-preserving  inclusion of finite von Neumann algebras $N\subseteq M$ and let $\mc A$ be an intermediate $\ast$-algebra (i.e.~$N\subseteq \mc A\subseteq M$) which is strongly dense in $M$. We will keep this setup fixed throughout the present  section and the symbols $N$, $\mc A$ and $M$ will therefore refer to objects specified above. Moreover, $\tau$ will denote the common trace on $N$ and $M$ and $H$ the GNS-space $L^2(M,\tau)$.
The aim of this section is to introduce a F{\o}lner type condition for $\mc A$, relative to the chain $N\subseteq \mc A\subseteq M$, and study its basic properties. 
As already mentioned, the von Neumann algebra $N$ will be thought of as a ``coefficient algebra'', and we advise the reader to keep the example $N=\mb C$ in mind in order to get a more intuitive picture. Recall from Section \ref{dim-theory-section}, that an $N$-submodule $\mc F$ in $M$ is called \emph{complete} if it is a Hilbert module with respect to the $N$-valued inner product arising from the trace-preserving conditional expectation $E\colon M\to N$.

\begin{Def}\label{strong-foelner-def}
The algebra $\mc A$ is said to have the strong F{\o}lner property with respect to  $N$ if for every finite set $T_1,\dots, T_r \in \mc A$ and every $\eps>0$ there exists a complete, non-zero, finitely generated projective $N$-submodule $\mc F\subseteq \mc A$ such that
\[
\frac{\dim_N(T_i^{-1}(\mc F)\cap \mc F)}{\dim_N(F)}>1-\eps \quad \text{ and } \quad  \|\varphi_{\mc F}-\tau  \|_{M_*}<\eps
\]
for all $i\in \{1,\dots, r\}$. Here $T_i^{-1}(\mc F)$ denotes the preimage of $\mc F$ under the left multiplication operator given by $T_i$. 
\end{Def}
This definition is an operator algebraic analogue of the F\o{}lner condition for discrete groups, where the almost invariant finite subset of the group is replaced by an almost invariant ``finite-dimensional'' $N$-submodule $\mc F$ in $\mc A$. In fact, putting $N=\CC$ one can easily check that the linear span of a subset $F$ in a group $\Gamma$ which is $\eps$-invariant under the action of another set $S$  gives rise to an almost invariant submodule $\mc F$ in the above sense (see Corollary \ref{grp-cor} for details) for any set of operators in $\CC\Gamma$ not supported outside of $S$. The condition regarding the trace approximation is trivially fulfilled in this case as $\varphi_{\mc F}(T)=\tau(T)$ for each $T\in L\Gamma$, a fact due to  $\CC\Gamma$ being spanned by a multiplicative set of orthogonal unitaries. Since this need not be the case for a general $\mc A$, we have to include the approximation property in order to compare the dimension of an $N$-submodule with the relative dimensions of its ``finite-dimensional approximations''; see Proposition \ref{matrix-approximation-prop} for the precise statement.

\begin{Rem}\label{left-right-rem}
Strictly speaking, the norm estimate in Definition \ref{strong-foelner-def} should read $\|\varphi_{\mc F}|_{M}-\tau\|_{M_*}<\eps$ but for notational convenience we will consistently suppress the restriction in the sequel. This should not lead to any confusion as the algebra on which the states are considered can always be read off the subscript on the norm. Moreover, the strong F{\o}lner property should, more precisely, be called the \emph{right} strong F{\o}lner property since we have chosen to use \emph{right} modules in the definition. However, if $\A$ has the right F{\o}lner property then it also has the corresponding left F{\o}lner property and vice versa. This can be seen by noting that if $T_1,\dots, T_r\in \A$ and $\eps>0$  are given and $\mc F\subseteq \mc A$ is a right $N$-module satisfying the conditions in Definition \ref{strong-foelner-def}, then $\mc F^*$ (the adjoint taken inside $\mc A$) is a complete, finitely generated projective \emph{left} Hilbert $N$-submodule in $M$ (the latter endowed with the inner product $\ip{x,y}_N=E(xy^*)$) which satisfies the conditions in the left version of Definition \ref{strong-foelner-def} for the set of operators $T_1^*,\dots, T_r^*$. Note, in particular, that this implies that the strong F{\o}lner property is stable under passing to opposite algebras; i.e.~if the tower $N\subseteq \mc A \subseteq M$ has the strong F{\o}lner property then the same is true for the tower $N^\op\subseteq \mc A^\op \subseteq M^\op$.

\end{Rem}
For our purposes the following slight reformulation of the F{\o}lner property will turn out convenient.

\begin{Prop}\label{strong-foelner-prop} 
The algebra $\A$ has the strong F{\o}lner property with respect to $N$ iff for any finite set $T_1,\dots,T_r\in \mc A$ 
 there exists a sequence of complete, non-zero, finitely generated, projective $N$-modules $\mc P_n\subseteq \mc A$  with submodules $\mc S_n\subseteq \mc P_n$ such that
\begin{enumerate}
\item[(i)] $T_i (\mc S_n) \subseteq \mc P_n$ for each $i\in\{1,\dots, r\}$ and each $n\in\mb N$;
\item[(ii)]  $\displaystyle \lim_{n\to \infty }\frac{\dim_N\mc S_n}{\dim_N\mc P_n}=1$;  
\item[(iii)] The sequence $\varphi_{\mc P_n}$ (restricted to $M$) converges uniformly to the trace $\tau$; i.e.~$\lim_{n\to \infty }\|\varphi_{\mc P_n}-\tau\|_{M_*}=0$.
\end{enumerate}
\end{Prop}
\begin{proof}
Clearly (i)-(iii) imply the F{\o}lner property. Conversely, if $\mc A$ has the F{\o}lner property and $T_1,\dots,T_r \in \mc A$ is given  then
for each $n\in \NN$ we get a complete, non-zero, finitely generated projective module $\mc F_n\subseteq \mc A $ such that 
\[
\frac{\dim_N(T_i^{-1}(\mc F_n)\cap \mc F_n)  }{\dim_N \mc F_n}>1-\frac{1}{rn} \quad \text{ and } \quad \|\varphi_{\mc F_n}-\tau\|_{M_*}<\frac{1}{n}.
\]
Putting $\mc P_n:=\mc F_n$ and $\mc S_n:=\cap_{i=1}^r T_i^{-1}(\mc F_n)\cap \mc F_n$ we clearly have (i) and (iii) satisfied and (ii) follows from Lemma \ref{intersection-lem}.
\end{proof}

\begin{Def}
If $\mc A$ has the strong F{\o}lner property and $\{T_1,\dots, T_r\}$ is a finite subset in $\mc A$, then we call a sequence $(\mc P_n,\mc S_n)$ with the properties 
in Proposition \ref{strong-foelner-prop} a \emph{strong F{\o}lner sequence} for the given set of operators. 
\end{Def}
The final result in this section shows that the von Neumann dimension can be approximated by relative dimensions in the presence of a strong F{\o}lner sequence.

\begin{Prop}\label{rel-to-vN-prop}
Let $\A$ have the strong F{\o}lner property relative to $N$ and let $(\mc P_n,\mc S_n)$ be a strong F{\o}lner sequence for an arbitrary finite set in $\mc A$. If $K\subseteq H^k$ is a closed subspace which is invariant under the diagonal right action of $M$ then $ \dim_{\mc P_n} K \underset{n\to \infty}{\To} \dim_M K$.
\end{Prop}
\begin{proof}
Denote by $P_K\in \MM_k(M)\subseteq B(H^k)$ the projection onto $K$. Fix an $n\in \NN$ and choose a basis $u_1,\dots, u_l\in \mc P_n$. Then the set
\[
\{ u_i \tens e_j \mid 1\leqslant i\leqslant l, 1\leqslant j\leqslant k \}
\]
is a basis for the amplification $\mc P_n^k=\mc P_n \tens \CC^k$; here $e_1,\dots, e_k$ denotes the standard basis in $\CC^k$. By computing the trace in this basis one easily gets
\[
\dim_{\mc P_n}K=\sum_{j=1}^k \varphi_{\mc P_n}((P_K)_{jj}) \underset{n\to\infty}{\To} \sum_{j=1}^k \tau((P_K)_{jj})=\dim_MK,
\]
where the convergence follows from Proposition \ref{strong-foelner-prop}.
\end{proof}

\section{Dimension flatness}
Throughout this section we consider again the  setup consisting of a trace-pre\-ser\-ving inclusion $N\subseteq M$ of tracial von Neumann algebras together with an intermediate $*$-algebra $N\subseteq \mc A \subseteq M$ which is weakly dense in $M$. We will also consider a $k$-fold amplification $H^k$ of the GNS-space $H:=L^2(M,\tau)$ and the natural left action of $\MM_k(\mc A)\subseteq \MM_k(M)$ thereon.  Our aim  is Theorem \ref{dim-flat-thm} which roughly says that when $\mc A$ is strongly F{\o}lner then the ring inclusion $\mc A \subseteq M$ is flat from the point of view of dimension theory. Before reaching our goal we need a few preparatory results.

\begin{Prop}\label{matrix-approximation-prop}
Assume that $\mc A$ has the strong F{\o}lner property relative to $N$ and let $T=(T_{ij})\in \MM_k(\mc A)$ be given. If $(\mc P_n,\mc S_n)$ is a strong F{\o}lner sequence for the set of matrix entries of $T$, then 
\[
\dim_M\ker(T)=\lim_{n\to \infty} \dim_{\mc P_n}\ker(T|_{\mc S_n^k}). 
\]
\end{Prop}
The proof is an extension of the corresponding argument in \cite{elek-zdc}.

\begin{proof}
By  L\"uck's dimension theorem \cite[Theorem 6.7]{Luck02} we have
\[
\dim_{N} \ker(T|_{\mc S_n^k}) + \dim_{N}  \rg(T|_{\mc S_n^k})=\dim_{N}\mc S_n^k=k \dim_N\mc S_n.
\]
Since $T_{ij}\mc S_n \subseteq \mc P_n$ we have ${\rg(T|_{\mc S_n^k})}\subseteq \mc P_n^k$ and $\ker(T|_{\mc S_n^k})\subseteq \mc P_n^k $ and by the basic properties of the relative dimension function (Proposition \ref{rel-dim-properties-amplified}) we now get
\begin{align}\label{estimate-0}
\dim_{\mc P_n} \ker(T|_{\mc S_n^k}) + \dim_{\mc P_n}{ \rg(T|_{\mc S_n^k}})=k \frac{\dim_N\mc S_n}{\dim_N \mc P_n}.
\end{align}
Denote by $P$ the kernel projection of $T\colon H^k\to H^k$ and by $R$ the projection onto the closure of its range. By Proposition \ref{rel-to-vN-prop}, for any $\varepsilon>0$ we can find $n_0\in {\mb N}$ such that for all $n\geqslant n_0$ we have 
\begin{align}
a_n&:=\dim_{\mc P_n}\ker(T|_{\mc S_n^k})\leqslant \dim_{\mc P_n}\ker(T) \leqslant \dim_M\ker(T)+\varepsilon; \label{estimate-1}  \\
b_n&:=\dim_{\mc P_n}{\rg(T|_{\mc S_n^k})} \leqslant \dim_{\mc P_n}\overline{\rg(T)}\leqslant \dim_M\overline{\rg(T)}+\varepsilon. \label{estimate-2}
\end{align}
By \eqref{estimate-0} we have $\lim_n(a_n+b_n)=k$ and by \cite[Theorem 1.12 (2)]{Luck02} 
\[
\dim_M\ker(T)+\dim_M\overline{\rg(T)}=k.
\]
 Our task is to prove that $a_n$ converges to $\dim_M\ker(T)$. If this were not the case, by passing to a subsequence we can assume that there exists $\varepsilon_0>0$ and an $n_1\in \NN$ such that $a_n\notin [\dim_M\ker(T)-\varepsilon_0,\dim_M\ker(T)+\varepsilon_0]$ for $n\geqslant n_1$. The estimates \eqref{estimate-1} and \eqref{estimate-2} imply the existence of an $n_2\in \NN $ such that 
 \begin{align*}
 a_n\leqslant \dim_M\ker(T)+\frac{\varepsilon_0}{2} \quad \text{ and } \quad b_n\leqslant \dim_M \overline{\rg(T)} + \frac{\varepsilon_0}{2}
\end{align*}
for $n\geqslant n_2$; hence we must have $a_n\leqslant \dim_M\ker(T)-\varepsilon_0$ for $n\geqslant\max\{n_1,n_2\}$. But then from this point on 
\[
a_n+b_n\leqslant \dim_M\ker(T)- \varepsilon_0 + \dim_M\overline{\rg(T)}+ \frac{\varepsilon_0}{2} = k-\frac{\varepsilon_0}{2},
\]
contradicting the fact that $\lim_n(a_n+b_n)=k$.
\end{proof}
Consider again the operator $T\in \MM_k(\mc A) \subseteq B(H^k)$ as well as its restriction $T_0\colon \mc A^k\to \mc A^k$.
\begin{Lemma}\label{closure-lem}
If $\mc A$ has the strong F{\o}lner property relative to $N$ then the closure of ${\ker(T_0)}$ in $H^k$ coincides with $\ker(T)$. 
\end{Lemma}
\begin{proof}
Let $(\mc P_n,\mc S_n)$ be a strong F{\o}lner sequence for the matrix coefficients of $T$. Denote by $P$ the kernel projection of $T$ and let $Q$ be the projection onto ${\overline{\ker(T_0)}}^\perp \cap \ker(T)$. We need to prove that $Q=0$. One easily checks that $\rg(Q)$ is a finitely generated, normal, right $M$-module
and it  is therefore enough to prove that $\dim_M\rg(Q)=0$. Denote by $R$ the projection onto the space $\overline{\ker(T_0)}$. Given  $\varepsilon>0$,  Proposition \ref{matrix-approximation-prop} provides an $n_0\in \NN$ such that for $n\geqslant n_0$ we have
\[
\dim_M\rg(P)=\dim_M\ker(T)\leqslant  \dim_{\mc P_n}\ker(T|_{\mc S_n^k})+\varepsilon \leqslant \dim_{\mc P_n}\rg(R)+\varepsilon,
\]
simply because $\ker(T|_{\mc S_n^k})\subseteq \ker(T_0)\subseteq \overline{\ker(T_0)}=\rg(R)$. Since $P=R+Q$, Proposition \ref{rel-to-vN-prop} implies that we eventually have
\begin{align*}
\dim_{\mc P_n}\rg(R)+\dim_{\mc P_n} \rg(Q) &=\dim_{\mc P_n}\rg(P)\\\
&\leqslant \dim_M\rg(P)+\varepsilon\\
& \leqslant  \dim_{\mc P_n}\rg(R)+2\varepsilon
\end{align*}
and hence $\dim_{\mc P_n}\rg(Q)\leqslant 2\varepsilon$ from a certain point on. Thus $\dim_M\rg(Q)=\lim_n \dim_{\mc P_n}\rg(Q)=0$.
\end{proof}
\begin{Rem}\label{left-right-rem-2}
If instead $T_0$ is given by \emph{right} multiplication with a matrix from $\MM_k(\mc A)$ then $T\in \diag(M)'=\MM_k(M')$  and hence commutes with the diagonal action  of $N$ from the left. By using the obvious variations of the above results for left modules we therefore also obtain $\overline{\ker(T_0)}=\ker(T)$ if $\mc A$ has the F{\o}lner property (see e.g.~Remark \ref{left-right-rem}). We will use this variant of the result in the proof of Theorem \ref{dim-flat-thm} which is formulated using dimension-theory for \emph{left} modules over $M$; this is done in order be consistent with the majority of the references (e.g.~\cite{Luck02, CS, sauer-betti-for-groupoids}) on $L^2$-Betti numbers in the homological algebraic context. 
\end{Rem}

We are now ready to state and prove the main theorem of this section. Recall, that $N\subseteq M$ is a trace-preserving inclusion of  finite von Neumann algebras and that $\mc A$ denotes an intermediate $*$-algebra which is weakly dense in $M$.
\begin{Thm}[Dimension flatness]\label{dim-flat-thm}
If $\mc A$ has the strong  F{\o}lner property relative to $N$ then the inclusion $\mc A \subseteq M$ is dimension flat; that is
\[
\dim_M \Tor_p^{\mc A} (M, X)=0
\]
for any $p\geqslant 1$ and any left $\mc A$-module $X$.
\end{Thm}
Note that if the ring inclusion $\mc A \subseteq M$ actually were flat (in the standard sense of homological algebra) then we would have $\Tor_p^{\mc A}(M,X)=0$ for every left $\mc A$-module $X$ and every $p\geqslant 1$. This need not be the case in our setup\footnote{For example, the tower $\CC\subseteq \CC[\ZZ\times \ZZ]\subseteq L(\ZZ\times \ZZ)$ has the strong F{\o}lner property (Corollary \ref{grp-cor}) but is not flat as $\Tor_1^{\CC[\ZZ\times \ZZ]}(L(\ZZ\times \ZZ),\CC)\neq \{0\}$. See e.g. \cite[Theorem 3.7]{LRS}.}, but from the point of view of the dimension function it looks as if it were the case --- hence the name ``dimension-flatness''. The first part of the proof of Theorem \ref{dim-flat-thm} consists of a reduction to the case when $X$ is finitely presented. This part is verbatim identical to the corresponding proof for groups due to L{\"u}ck (see \cite[Theorem 6.37]{Luck02}), but we include it here for the sake of completeness. 

\begin{proof}
Let an arbitrary left $\mc A$-module $X$ be given and choose a short exact sequence $0\to Y \to F \to X \to 0$ of $\mc A$-modules in which $F$ is  free. Then the corresponding long exact $\Tor$-sequence gives
\[
\Tor_{p+1}^{\mc A}(M, X)\simeq \Tor_{p}^{\mc A}(M,Y),
\]
and hence it suffices to prove that $\dim_M \Tor_1^{\mc A} (M, X)=0$ for arbitrary $X$. Recall that $\Tor$ commutes with direct limits and that the dimension function $\dim_M(-)$ is also well behaved with respect to direct limits \cite[Theorem 6.13]{Luck02}); seeing that an arbitrary module is the directed union of its finitely generated submodules we may therefore assume $X$ to be finitely generated. Hence we can find a short exact sequence $0\to Y\to F\to X\to 0$ with $F$ finitely generated and free. The module $Y$ is the directed union of its system of finitely generated submodules $(Y_j)_{j\in J}$ and therefore $X$ can be realized as the direct limit $\varinjlim_j F/Y_j$. Since each of the modules $F/Y_j$ is finitely presented by construction this shows that it suffices to treat the case where $X$ is a finitely presented module. In this case we may therefore choose a presentation of the form
\[
\mc A^k \overset{T_0}{\To} \mc A^l \To X\To 0.
\] 
This presentation can be continued to a free resolution 
\[
\cdots\overset{S_0}{\To}\mc A^k \overset{T_0}{\To} \mc A^l \To X\To 0
\] 
of $X$ that can be used to compute the $\Tor$ functor;  in particular we get 
\[
\Tor_1^{\mc A}(M, X)\simeq \frac{\ker(\id_M\tens_{\mc A} T_0)}{\rg(\id_M\tens_{\mc A} S_0)}.
\]
Denote by $T_0^\vN\colon M^k\to M^l$ the map induced by $\id_M\tens_{\mc A} T_0$ under the natural identification $M\odot_{\mc A} \mc A^i=M^i$ ($i=k,l$) and 
by $T\colon H^k \to H^l$ its continuous extension\footnote{Recall that ``$\odot$'' denotes the algebraic tensor product and $H$ the GNS-space $L^2(M,\tau)$.}.  
Since $M\odot_{\mc A} -$ is right exact and $S_0$ surjects onto $\ker(T_0)$, 
we see that  $\rg(\id_M\tens_{\mc A} S_0)\subseteq M\odot_{\mc A} \mc A^k$ is identified with the $M$-submodule $M\ker(T_0)$ in $M^k$ generated by $\ker(T_0)$;  thus
\begin{align}\label{dim-Tor}
\dim_M \Tor_1^{\mc A}(M, X)=\dim_M \ker(T_0^\vN) - \dim_M (M\ker(T_0)).
\end{align}
We now claim that
\begin{align}\label{closure-eqs}
\dim_M \ker T_0^\vN =\dim_M \ker(T)  \quad \text{ and } \quad \dim_M(M\ker(T_0))=\dim_M \overline{\ker(T_0)},
\end{align}
where the closure is taken in the Hilbert space $H^k$. The first equality follows from the fact that the $L^2$-completion functor has an exact and dimension preserving inverse \cite[Theorem 6.24]{Luck02}. 
By Lemma \ref{same-closure-lem}, the dimension of a submodule in $M^k$ coincides with the dimension of its closure in $H^k$, so to prove the  the second equality it suffices to prove that $M\ker(T_0)$ and $\ker(T_0)$ have the same closure in $H^k$. But this  follows from the following simple approximation argument:
Take $x\in M$ and $a\in \ker(T_0)$ and consider $xa\in M\ker(T_0)$. By picking a net $x_\alpha\in \mc A$ converging in the strong operator topology to $x$ we obtain that $\|xa-x_\alpha a\|_2\to 0$ and $x_\alpha a\in \ker(T_0)$; hence $M\ker(T_0)\subseteq \overline{\ker(T_0)}$ and the proof of \eqref{closure-eqs} is complete. By \eqref{dim-Tor} and \eqref{closure-eqs} it is therefore sufficient to prove that $\overline{\ker(T_0)}=\ker(T)$. The adjoint operator $T^*\colon H^l \to H^k$ is the extension of the operator  $T_0^*\colon\mc A^l \to \mc A^k$ which multiplies from the right with the adjoint of the matrix defining $T_0$. From this we obtain $\ker(T_0^*T_0^{\phantom{*}})=\ker(T_0)$ and since $\ker(T^*T)=\ker(T)$ it is equivalent to prove that
\[
\overline{\ker(T_0^*T_0^{\phantom{*}})}=\ker(T^* T),
\]
but this follows from Lemma \ref{closure-lem} and Remark \ref{left-right-rem-2}.
\end{proof}
\begin{Rem}
A careful examination of the results in this section reveals that the dimension-flatness theorem can actually be obtained under slightly less restrictive assumptions. Namely, the proof goes through as long as $N\subseteq \mc A \subseteq M$ satisfies the requirements (i)-(iii) from Proposition \ref{strong-foelner-prop}, but with the uniform convergence in (iii) replaced by weak$^*$ convergence. 
However, for the results in the following section the uniform convergence will be of importance which is the reason for it being included in the definition of the strong F{\o}lner property.
\end{Rem}

\section{Operator algebraic amenability}
In this section we explore the connection between the strong F{\o}lner property and the existing operator algebraic notions of amenability. Consider therefore again a trace-preserving inclusion $N\subseteq M$ of finite von Neumann algebras and an intermediate $*$-algebra $\mc A$ which is weakly dense in $M$. The main goal  is to prove the following result:
\begin{Thm}\label{injectivity-thm}
If  $\mc A$ has the strong F{\o}lner property relative to $N$ then $M$ is amena\-ble relative to $N$. In particular, $M$ is injective if  $\mc A$ has the strong F{\o}lner property and $N$ is injective.
\end{Thm}
The notion of relative amenability for von Neumann algebras dates back to Popa's  work in \cite{popa-correspondences}; we briefly recall the basics here following the exposition in \cite{ozawa-popa-on-a-class}. Consider again the trace-preserving inclusion $N\subseteq M$ of finite von Neumann algebras and denote by $Q:=\langle M,e_N\rangle = (JNJ)'\cap B(L^2(M))$ the basic construction. Recall that $Q$  has a unique normal, semifinite tracial weight $\Psi\colon Q_+\to [0,\infty]$ with the property that
\[
\Psi(ae_N b)=\tau(ab) \text{ for all } a,b\in N.
\]
One way to construct the trace $\Psi$ is as follows: Since any normal representation of $N^\op$ is an amplification followed by a reduction, there exists a projection $q\in B(\ell^2(\NN)\otimes L^2(N)) $ and  a right $N$-equivariant unitary identification 
\[
L^2(M)=q(\ell^2(\NN)\otimes L^2(N)),
\]
where the right hand side is endowed with the diagonal right $N$-action. This induces an identification $Q:=(JNJ)'\cap B(L^2(M))=q (B(\ell^2(\NN)\otimes N  ))q$ and the trace $\Psi$ is simply the pull back of the restriction of $\Tr\otimes \tau$ under this identification. 

\begin{Thm}[{\cite[Theorem 2.1]{ozawa-popa-on-a-class}}]\label{popa-ozawa-thm}
The following conditions are equivalent:
\begin{itemize} 
\item[(i)] There exists a state $\varphi\colon Q\to \CC$ such that $\varphi |_M=\tau$ and  $\varphi(xT)=\varphi(Tx)$ for all $x\in M$ and $T\in Q$.
\item[(ii)] There exists a conditional expectation $E\colon Q\to M$.
\item[(iii)] There exists a net $\xi_n \in L^2(Q,\Psi)$ such that $\lim_n \ip{\xi_n,x\xi_n}_{L^2(Q,\Psi)}=\tau(x)$ and $\lim_n\|[x,\xi_n]\|_{2,\Psi}=0$ for every $x\in M$.
\end{itemize}
If $N\subseteq M$ satisfies one, and hence all, of these conditions then $M$ is  said to be amenable relative to $N$ (or $N$ to be coamenable in $M$). 
\end{Thm}

Note that (ii) in the above theorem gives that amenability of $M$ relative to $\CC$ is equivalent to amenability (a.k.a.~injectivity) of $M$. Note also that if $\mc F\subseteq M$ is a complete, finitely generated, projective right $N$-module then the projection $P_F$ onto its closure $F$ in $L^2(M)$ is an element in $Q$ and
\[
\dim_N(\mc F)=\Psi(P_F).
\]
This follows from the construction of $\Psi$, the equality $\dim_N(\mc F)=\dim_N(F)$ and the definition of the von Neumann dimension for normal $N$-modules; see \cite[Definition 1.8]{Luck02} and the comments following it. More generally, for any operator $T\in Q$ we have $\Psi(P_FTP_F)=\Tr_N(P_FTP_F)$ where $\Tr_N$ is the trace on $\End_N(F)$ considered in Section \ref{dim-theory-section} (see also Lemma \ref{trace-lemma}).
We are now ready to give the proof of Theorem \ref{injectivity-thm}

\begin{proof}[Proof of Theorem \ref{injectivity-thm}.]
Assume that $N\subseteq \mc A\subseteq M$ has the F{\o}lner property.  Since $M_*$ is assumed separable the unit ball $(M)_1$ is separable and metrizable
 for the strong operator topology and by Kaplansky's density theorem the unit ball $(\mc  A)_1$ is strongly dense in $(M)_1$. We may therefore choose a sequence $\{T_i\}_{i=1}^\infty$ in $(\mc A)_1$ which is strongly dense in $(M)_1$ and upon adding further operators to this sequence we may assume that $\{T_i\}_{i=1}^\infty$ is $*$-stable.
Since $\mc A$ satisfies the strong F{\o}lner condition we can, for each  $T_1,\dots, T_n$, find a complete, finitely generated, projective $N$-submodule $\mc P_n \subseteq \mc A$ and a submodule $\mc S_n \subseteq \mc P_n$ such that
\begin{itemize}
\item   $T_i(\mc S_n)\subseteq \mc P_n$ for all $i\in \{1,\dots, n\}$,
\item $\frac{\dim_N \mc S_n}{\dim_N\mc P_n}\geqslant 1-\frac{1}{n}$,
\item$ \|\varphi_{\mc P_n}- \tau \|_{M_*}\leqslant \frac{1}{n}$.
\end{itemize} 
In what follows we denote by $P_n\in B(H)$ the orthogonal projection onto the closure (in $H$) of $\mc P_n$, by $S_n$ the orthogonal projection onto the closure of $\mc S_n$ and by $S_n^{\perp}$ the difference $P_n-S_n$. Since $\mc P_n$ is complete and finitely generated projective, the discussion preceding the proof implies  that $\Psi(P_n)=\dim_N\mc P_N$; in particular $P_n\in L^2(Q,\Psi)$. We aim at proving that the unit vectors
\[
\xi_n :=\frac{1}{\sqrt{\dim_N\mc P_n}}P_n \in L^2(Q,\Psi)
\]
satisfy condition (iii) of Theorem \ref{popa-ozawa-thm}. The trace approximation is automatic since
\[
\ip{\xi_n,x\xi_n}_{L^2(Q,\Psi)}=\frac{\Psi(P_nx P_n)}{\dim_N\mc P_n}=\frac{\Tr_N(P_nx P_n)}{\dim_N\mc P_n}=\varphi_{\mc P_n}(x)\underset{n\to\infty}{\To} \tau(x)
\]
for any $x\in M$. To prove the asymptotic commutation property, consider first $x=T_{i_1}\in \{T_i\}_{i=1}^\infty$. Then $x^*=T_{i_2}$ for some $i_2\in \NN$
so for for $n\geqslant \max\{i_1,i_2\}$ we have $P_nxS_n=xS_n$ and $P_nx^*S_n=x^*S_n$. Hence
\begin{align*}
\|xP_n - P_n x\|_{2,\Psi}^2 &=\Psi((xP_n-P_nx)^*(xP_n-P_nx))\\
&= \Psi(P_nx^*xP_n-P_nx^*P_n x -x^*P_nxP_n + x^*P_nx  ) \\
&=\Psi(P_nx^*x P_n- P_n x^* P_n xP_n -P_nxP_n x^*P_n +  P_nxx^*P_n)\\
&=\Psi\left((P_nx^*x P_n- P_n x^* P_n xP_n -P_nxP_n x^*P_n +  P_nxx^*P_n)S_n^\perp\right)\\
&=\Psi\left(S_n^\perp(P_nx^*x P_n- P_n x^* P_n xP_n -P_nxP_n x^*P_n +  P_nxx^*P_n)S_n^\perp\right)\\
&=\left|  \ip{S_n^\perp, (P_nx^*x P_n- P_n x^* P_n xP_n -P_nxP_n x^*P_n +  P_nxx^*P_n)S_n^\perp }_{L^2(Q,\Psi)} \right|\\
&\leqslant \|S_n^\perp\|_{2,\Psi}\|(P_nx^*x P_n- P_n x^* P_n xP_n -P_nxP_n x^*P_n +  P_nxx^*P_n)S_n^\perp\|_{2,\Psi}\\
&\leqslant \|S_n^\perp\|_{2,\Psi}^2\|P_nx^*x P_n- P_n x^* P_n xP_n -P_nxP_n x^*P_n +  P_nxx^*P_n\|_{\infty}\\
&\leqslant 4 \Psi(S_n^{\perp})\\
&= 4(\dim_N\mc P_n-\dim_N\mc S_n).
\end{align*}
Hence
\[
\|[x,\xi_n]\|_{2,\Psi}=\frac{\|xP_n-P_nx\|_{2,\Psi}}{\sqrt{\dim_N \mc P_n}}\leqslant 2\sqrt{\frac{\dim_N \mc P_n-\dim_N\mc S_n}{\dim_N \mc S_n}}\underset{n\to\infty}{\To} 0.
\]
Now the general case follows from this and an approximation argument: let $x\in M$ and $\varps>0$ be given and assume without loss of generality that $\|x\|_\infty\leqslant 1$. Choose a net $x_\alpha \in \{T_i\}_{i=1}^\infty$ converging strongly to $x$. Then $x_\alpha$ also converges to $x$ in $L^2(M,\tau)$ and we may therefore choose an $\alpha$ such that 
\[
\|x-x_\alpha\|_{2,\tau}<\eps.
\]
We now have
\[
\|x\xi_n-\xi_n x\|_{2,\Psi}\leqslant \|(x-x_\alpha)\xi_n \|_{2,\Psi} + \|x_\alpha \xi_n -\xi_n x_\alpha \|_{2,\Psi} +\| \xi_n(x_\alpha -x) \|_{2,\Psi}
\]
Considering the first term we get
\begin{align*}
\|(x-x_\alpha)\xi_n \|_{2,\Psi}^2 &= \frac{\Psi(P_n (x-x_\alpha)^*(x-x_\alpha)P_n   )}{\dim_N\mc P_n}\\
&=\varphi_{\mc P_n}((x-x_\alpha)^*(x-x_\alpha))\\
&\leqslant |(\varphi_{\mc P_n}-\tau)((x-x_\alpha)^*(x-x_\alpha)) | +\tau((x-x_\alpha)^*(x-x_\alpha)) \\
&\leqslant \|\varphi_{\mc P_n}-\tau \|_{M_*} \|x-x_{\alpha}\|_\infty^2 +\|x-x_\alpha\|_{2,\tau}^2\\
&\leqslant \frac{2}{n} +\eps^2.
\end{align*}
Thus, for a suitably chosen $n_1\in \NN$ we have $\|(x-x_\alpha)\xi_n \|_{2,\Psi}^2\leqslant 2\varps$ for all $n\geqslant n_1$. Considering the third term we get, in a completely similar manner, an $n_3\in \NN$ such that $\| \xi_n(x_\alpha -x) \|_{2,\Psi}\leqslant 2\eps$ for $n\geqslant n_3$. Since $x_\alpha \in \{T_i\}_{i=1}^\infty$, the second term converges to zero by what was shown in the first part of the proof, and hence there exists an $n_2\in \NN$ such that $\|x_\alpha \xi_n -\xi_n x_\alpha \|_{2,\Psi}\leqslant \eps$ for $n\geqslant n_2$. Thus, for $n\geqslant \max\{n_1,n_2,n_3\}$ we have $\|x\xi_n-\xi_n x \|_{2,\Psi}\leqslant 5\eps$ and since $\eps>0$ was arbitrary this completes the proof of the amenability of $M$ relative to $N$.\\

We now just have to prove the final claim regarding the injectivity of $M$. So assume that $N$ is injective and, as before, that $\mc A\subseteq M$ satisfies the F{\o}lner condition relative to $N$. Injectivity of $N$ simply means that $N$ is amenable relative to the subalgebra $\CC\subseteq N $ and by what was just proven $M$ is amenable relative to $N$ as well. By  \cite[Theorem 3.2.4]{popa-correspondences} this means that $M$ is also amenable relative to $\CC$; i.e~$M$ is injective.

\end{proof}

\section{\texorpdfstring{Examples and applications to $L^2$-Betti numbers}{Examples and applications to L2-Betti numbers}}\label{example-section}
The first goal in this section is to show how our techniques can be used to unify the dimension-flatness results known for amenable groups and quantum groups and, more generally, to provide a vanishing result in the context of operator algebraic $L^2$-Betti numbers. 
Secondly, in order to convince the reader that the class of algebras with the strong F{\o}lner property is rich and extends beyond the ones arising from amenable groups an quantum groups, we give a number of such examples and show how they yield results about vanishing of $L^2$-Betti numbers in various situations.   As a first application we re-obtain the original dimension flatness result of L{\"u}ck \cite[Theorem 6.37]{Luck02}. 

\begin{Cor}[$L^2$-Betti numbers of amenable groups]\label{grp-cor}
If $\Gamma$ is a countable, discrete, amenable group then  $\CC\Gamma$ has the strong F{\o}lner property with respect to the tower $\CC\subseteq \CC\Gamma \subseteq L\Gamma$. In particular, the inclusion $\CC\Gamma \subseteq L\Gamma$ is dimension flat and  $\beta_p^{(2)}(\Gamma)=0$ for $p\geqslant 1$.  
\end{Cor}
\begin{proof}
Let $T_1,\dots, T_r\in \CC\Gamma$  be given and put $S=\bigcup_{i=1}^r \supp(T_i)\subseteq \Gamma$. Since $\Gamma$ is amenable,  we can find (see e.g.~\cite[F.6.8]{benedetti}) a sequence of finite subset $F_n\subseteq \Gamma$ such that
\[
\frac{\card(\partial_S(F_n))}{\card(F_n)}\underset{n\to \infty}{\To} 0, 
\]
where $
\partial_S(F):=\{ \gamma \in F \mid \exists s\in S: s\gamma \notin F \}.$
Then putting $\mc P_n=\CC [F_n]\subseteq \CC\Gamma$ and $\mc S_n=\CC[F_n\setminus \partial_S(F_n)]$ we clearly have $T_i(\mc S_n)\subseteq \mc P_n$ for all $i\in \{1,\dots, r\}$ and 
\[
\frac{\dim_{\CC} \mc S_n  }{\dim_\CC \mc P_n}= \frac{ \card(F_n) - \card(\partial_S(F))  }{\card(F_n)}\underset{n\to \infty}{\To} 1.
\]
Denoting by $\rho$ the right regular representation of $\Gamma$ we get for $T\in L\Gamma$ 

\begin{align*}
\varphi_{\mc F}(T)=\frac{\sum_{\gamma\in F} \ip{\delta_\gamma, T\delta_\gamma}}{\card(F)}=\frac{\sum_{\gamma\in F} \ip{\delta_e, \rho_{\gamma}^*T\rho_\gamma\delta_e}}{\card(F)} = \frac{\sum_{\gamma\in F} \ip{\delta_e,T\delta_e}}{\card(F)}=\tau(T).
\end{align*}
Hence $\CC\Gamma$ has the strong F{\o}lner property and dimension flatness therefore follows from Theorem \ref{dim-flat-thm}. The fact that the $L^2$-Betti numbers vanish in positive degrees follows from this since
\[
\beta_p^{(2)}(\Gamma)=\dim_{L\Gamma}\Tor_p^{\CC\Gamma}(L\Gamma,\CC)=0. \qedhere
\]
\end{proof}
Next we turn our attention to the class of compact/discrete quantum groups. Since this is a bit of an aside we shall not elaborate on the basics of quantum group theory; for back ground material the reader is referred to the detailed survey articles \cite{wor-cp-qgrps}, \cite{kustermans-tuset} and \cite{van-daele}.  In what follows, we denote by $\widehat{\GG}$ a discrete quantum group of Kac type and by $\GG=(C(\GG),\Delta_\GG)$ its compact dual. Associated with $\GG$ is a Hopf $*$-algebra $\Pol(\GG)$ which is naturally included into a finite von Neumann algebra $L^\infty(\GG)$. Furthermore, the von Neumann algebra $L^\infty(\GG)$ carries a natural trace $\tau$, namely the Haar state associated with $\GG$. 
In \cite{coamenable-betti} and \cite{quantum-betti}  $L^2$-Betti numbers were  studied in the context of discrete quantum groups of Kac type and we re-obtain \cite[Theorem 6.1]{coamenable-betti} by setting $N=\CC$, $\mc A= \Pol(\GG)$ and $M=L^\infty(\GG)$: 

\begin{Cor}[$L^2$-Betti numbers of amenable quantum groups]\label{qgrp-cor}
If $\widehat{\GG}$ is a dis\-cre\-te amenable quantum group of Kac type and $\GG$ is its compact dual then $\Pol(\GG)$ has the strong F{\o}lner property with respect to the tower $\CC\subseteq \Pol(\GG)\subseteq L^\infty(\GG)$. In particular, the inclusion $\Pol(\GG)\subseteq L^\infty (\GG)$ is dimension flat and $\beta_p^{(2)}(\widehat{\GG})=0$ for $p\geqslant 1$.
\end{Cor}
\begin{proof}
Choose a complete set $(u^\alpha)_{\alpha\in I}$ of representatives for the set of equivalence classes of irreducible unitary corepresentations of $\GG$ and denote by $n_\alpha$ the matrix size of $u^\alpha$. Recall that the corresponding matrix coefficients $u_{ij}^\alpha$ constitute a linear basis for $\Pol(\GG)$ and that the normalized matrix coefficients $\sqrt{n_\alpha}u_{ij}^\alpha$ constitute an orthonormal basis in $L^2(\GG):=L^2(\Pol(\GG),\tau)$. Let $T_1,\dots, T_r\in \Pol(\GG)$ and $\eps>0$ be given and denote by $S$ the joint support of $T_1,\dots, T_r$; that is 
\[
S=\{ u^\gamma \mid \exists l\in \{1,\dots, r \} \exists p,q\in \{1,\dots, n_\gamma\} :  \tau(u_{pq}^{\gamma *}T_l)\neq 0  \}.
\]
In other words, $S$ consists of all the irreducible corepresentations that has at least one of its matrix coefficients entering in the linear expansion of one of the operators in question. According to the quantum F{\o}lner condition \cite[Corollary 4.10]{coamenable-betti} we can find a sequence of finite subset $F_k\subseteq (u^\alpha)_{\alpha\in I}$ such that
\[
\sum_{u^\alpha \in \partial_S(F_k)}n_\alpha^2 < \frac{1}{k} \sum_{u^\alpha \in F_k}n_\alpha^2,
\]
where the boundary $\partial_S(F)$ is as in \cite[Definition 3.2]{coamenable-betti}. Define
\begin{align*}
\mc P_k&:=\spann_\CC\{u_{ij}^\alpha \mid u^\alpha\in F, 1\leqslant i,j\leqslant n_\alpha\}; \\ 
 \mc S_k&:=\spann_\CC\{u_{ij}^\alpha \mid u^\alpha \in F\setminus \partial_S(F), 1\leqslant i,j\leqslant n_\alpha\} .
\end{align*}
By construction we now have $T_l(\mc S_k)\subseteq \mc P_k$ for every $l\in \{1,\dots, r\}$ and that
\[
\frac{\dim_\CC \mc S_k}{\dim_\CC \mc P_k}= \frac{\sum_{u^\alpha \in F_k\setminus \partial_S(F_k)} n_\alpha^2 }{\sum_{\alpha\in F_k}n_\alpha^2} \underset{n\to\infty}{\To} 1.
\]
Furthermore, denoting by $P_{\mc P_n}\in B(L^2(\GG))$ the projection onto the closed subspace $\mc P_n\subseteq L^2(\GG)$ we get
\begin{align*}
\Tr_\CC(P_{\mc P_n} T P_{\mc P_n}  ) &= \sum_{\alpha \in F}\sum_{i,j=1}^{n_\alpha}\ip{\sqrt{n_\alpha}u_{ij}^\alpha, T\sqrt{n_\alpha}u_{ij}^\alpha}\\
&= \sum_{\alpha \in F}\sum_{i,j=1}^{n_\alpha} n_\alpha\tau(u_{ij}^{\alpha*}Tu_{ij}^\alpha)\\
&=\sum_{\alpha \in F} n_\alpha \tau\left(\left(\sum_{i,j=1}^{n_{\alpha}}u_{ij}^{\alpha}u_{ij}^{\alpha*}\right)T\right)\\
&=\tau(T)\sum_{\alpha \in F}n_{\alpha}^2\\
&=\tau(T)\dim_{\CC}\mc P_n,
\end{align*}
for every $T\in L^\infty(\GG)$;  hence $\varphi_{\mc P_n}(T)=\tau(T)$. Thus $\Pol(\GG)$ has the strong F{\o}lner property and dimension flatness follows from Theorem \ref{dim-flat-thm}. The vanishing of the $L^2$-Betti numbers follows from dimension flatness  since 
\[
\beta_p^{(2)}(\widehat{\GG}):=\dim_{L^\infty(\GG)}\Tor_p^{\Pol(\GG)}(L^\infty(\GG),\CC)=0. \qedhere
\] 
\end{proof}

Note that Corollary \ref{grp-cor} is actually just a special case of Corollary \ref{qgrp-cor}, but since we expect most readers to be more familiar with groups than quantum groups we singled out this case in form of Corollary \ref{grp-cor} and its proof. \\
\begin{Rem}
From Theorem \ref{injectivity-thm} and Corollary \ref{grp-cor} we re-obtain the classical fact that the group von Neumann algebra $L\Gamma$ is hyperfinite when $\Gamma$ is an amenable discrete group. Furthermore, we see that the group algebra $\CC\Gamma$ associated with a discrete group $\Gamma$ has the strong F{\o}lner property (relative to the chain $\CC\subseteq \CC\Gamma \subseteq L\Gamma$) if and only if $\Gamma$ is amenable: We already saw in Corollary \ref{grp-cor} that amenability implies the strong F{\o}lner property. Conversely, if $\CC\Gamma$ has the strong F{\o}lner property then $L\Gamma$ is hyperfinite by Theorem \ref{injectivity-thm} and thus $\Gamma$ is amenable (see e.g.~\cite[Theorem 2.6.8 and 9.3.3]{brown-ozawa}). By a similar argument, using \cite[Theorem 4.5]{ruan-amenability} and Corollary \ref{qgrp-cor}, we obtain that a discrete quantum group $\widehat{\GG}$ of Kac type is amenable if and only if the dual Hopf $*$-algebra $\Pol(\GG)$ has the strong F{\o}lner property relative to the chain $\CC\subseteq \Pol(\GG) \subseteq L^\infty(\GG)$. Thus, for algebras arising from groups and quantum groups the strong F{\o}lner condition coincides exactly with the notion of amenability of the underlying object.
\end{Rem}

Next we take a step up in generality by considering $L^2$-Betti numbers in a purely operator algebraic context.
In \cite{CS}, Connes and Shlyakhtenko introduce a notion of $L^2$-homology and $L^2$-Betti numbers for a dense $*$-subalgebra in a tracial von Neumann algebra $(M,\tau)$. More precisely, if $\mc A \subseteq M$ is a weakly dense unital $*$-subalgebra its $L^2$-homology and $L^2$-Betti numbers are defined,  respectively, as
\[
H_p^{(2)}(\mc A)=\Tor_p^{\mc A \odot \mc A^\op}(M\bar{\tens}M^{\op}, \mc A ) \ \text{ and } \ \beta_p^{(2)}(\mc A,\tau)=\dim_{M\bar{\tens}M^{\op}}H^{(2)}_p(\mc A).
\]
This definition extends the definition for groups \cite[Proposition 2.3]{CS} by means of the formula $\beta_p^{(2)}(\CC\Gamma,\tau)=\beta_p^{(2)}(\Gamma)$ and it also fits with the notion for quantum groups studied above  \cite[Theorem 4.1]{quantum-betti}.  In order to apply our techniques in the Connes-Shlyakhtenko setting we first need to prove that the strong F{\o}lner property is preserved under forming algebraic tensor products.

\begin{Prop}\label{tensor-foelner}
Let $\mc A$ and $\mc B$ be  dense $*$-subalgebras in the tracial von Neumann algebras $(M,\tau)$ and $(N,\rho)$, respectively.  
If $\mc A$ and $\mc B$ have the strong F{\o}lner property with respect to the towers $\CC\subseteq \mc A \subseteq M$ and $\CC\subseteq \mc B \subseteq N$  then so does $\mc A \odot \mc B$ with respect to the tower $\CC\subseteq \mc A \odot \mc B \subseteq M\bar{\tens} N$.
\end{Prop}
\begin{proof}
Let $T_1,\dots, T_r \in \mc A \odot \mc B$ be given and write each $T_i$ as $\sum_{k=1}^{N_i} a_k^{(i)}\tens b_{k}^{(i)}$. Choose a strong F{\o}lner sequence $\mc P_n'\subseteq \mc A$, with associated subspaces $\mc S_n'$, for the $a_k^{(i)}$'s and a strong F{\o}lner sequence  $\mc P_n''\subseteq \mc B$, with associated subspaces $\mc S_n''$,  for the $b_k^{(i)}$'s . Putting $\mc P_n:= \mc P_n'\odot \mc P_n'' $ and $\mc S_n:=\mc S_n'\otimes \mc S_n'' $ we clearly
have that the sequence $(\mc P_n, \mc S_n)$ satisfies the first to requirements in Proposition \ref{strong-foelner-prop}, and since $\varphi_{\mc P_n}=\varphi_{\mc P_n'}\tens \varphi_{\mc P_n''}$ we get
\begin{align*}
\|\varphi_{\mc P_n} -\tau\tens \rho\|_{(M\bar{\tens} N)_*} &= \|\varphi_{\mc P_n'}\tens \varphi_{\mc P_n''} - \tau\tens \varphi_{\mc P_n''} + \tau\tens\varphi_{\mc P_n''} -\tau\tens \rho \|_{(M\bar{\tens} N)_*}   \\
&\leqslant \| \varphi_{\mc P_n'}- \tau \|_{M_*} + \|\rho-\varphi_{\mc P_n''}\|_{N_*} \underset{n\to \infty}{\To} 0.
\end{align*}
Thus, $(\mc P_n, \mc S_n)$ is a strong F{\o}lner sequence for $T_1,\dots, T_r$.
\end{proof}

\begin{Cor}\label{CS-vanishing}
If $\mc A$ has the strong F{\o}lner property with respect to the tower $\CC\subseteq \mc A \subseteq M$ then the Connes-Shlyakhtenko $L^2$-Betti numbers of $(\mc A,\tau)$ vanish in positive degrees.
\end{Cor}
\begin{proof}
Since $\mc A$ is strongly F{\o}lner so is the opposite algebra $\mc A^\op$ (see e.g.~Remark \ref{left-right-rem}) and by Proposition \ref{tensor-foelner} the same is true for $\mc A \odot \mc A^\op$. Hence by Theorem \ref{dim-flat-thm} the inclusion $\mc A\odot \mc A^\op \subseteq M\bar{\otimes} M^\op$ is dimension flat so in particular
\[
\bet_p(\mc A, \tau)= \dim_{M\bar{\otimes} M^\op}\Tor_p^{\mc A \otimes \mc A^\op}(M\bar{\otimes} M^\op, \mc A)=0
\]
for all $p\geqslant 1$.
\end{proof}

\begin{Cor}[Actions of amenable groups on probability spaces]\label{action-cor}
Let  $\Gamma$ be a discrete, countable amenable group acting on a probability space $(X,\mu)$ in a probability measure-pre\-ser\-ving way and denote by $\alpha$ the induced action $ \Gamma \to \Aut(L^\infty(X))$. Consider the crossed product von Neumann algebra $M=L^\infty(X)\rtimes \Gamma$ and its subalgebras $N=L^\infty(X)$ and $\mc A = L^\infty(X)\rtimes_{\mathrm{alg}}\Gamma$. Then $\mc A$ has the strong F{\o}lner property relative to $N$ and hence the inclusion $L^\infty(X)\rtimes_{\mathrm{alg}}\Gamma \subseteq L^\infty(X)\rtimes\Gamma$ is dimension-flat. Moreover, if the action of $\Gamma$ is essentially free, then the $L^2$-Betti numbers (in the sense of Sauer \cite{sauer-betti-for-groupoids}) of the groupoid defined by the action vanish in positive degrees.
\end{Cor}
We denote by $(u_g)_{g\in \Gamma}$ the natural unitaries on $L^2(X)$ implementing the action; i.e.~$\alpha_g(a)=u_g^{\phantom{*}}au_g^*$ for all $g\in \Gamma$ and all $a\in L^\infty(X)$.
\begin{proof}
For any finite set $F\subseteq \Gamma$, the right $L^\infty(X)$-module
\beqn
\mc F = \left.\left\{\sum_{g\in F} a_g u_g\,\right|\, a_g \in L^\infty(X)\right\}
\eeqn
is free of dimension $\card(F)$ and the inner product on $\mc F$ arising from the conditional expectation $E\colon M\to N$ is given by
\[
\ip{\sum_{g\in F}a_gu_g, \sum_{g\in F} b_g u_g}_N=\sum_{g\in F} \alpha_{g^{-1}}(a_g^*b_g^{\phantom{*}}).
\]
From this it follows that the map $(L^\infty(X)^{\card(F)},\ip{\cdot, \cdot}_\st)\to (\mc F, \ip{\cdot, \cdot}_N)$ given by $(a_g)_{g\in F}\mapsto \sum_{g\in F}a_g u_g$ is a unitary isomorphism and that $u_1,\dots, u_n$ is a basis for $\mc F$. In particular, $\mc F$ is complete and finitely generated projective. The trace-approximation condition is automatic, because for $T\in M$ we have
\[
\varphi_{\mc F}(T)=\frac{\sum_{g\in F} \ip{u_g, Tu_g} }{\dim_N \mc F}= \frac{\sum_{g\in F} \tau(u_g^*Tu_g^{\phantom{*}}) }{\card(F)}=\tau(T).
\]
Hence the strong F{\o}lner condition for $\mc A$ relative to $N$ follows from the classical F{\o}lner condition for $\Gamma$ and the inclusion $L^\infty(X)\rtimes_{\mathrm{alg}}\Gamma \subseteq L^\infty(X)\rtimes\Gamma$ is therefore dimension-flat. In the case of an essentially free action we combine this with \cite[Lemma 5.4]{sauer-betti-for-groupoids} and get that the $L^2$-Betti numbers of the groupoid defined by an action of an amenable group vanish.
\end{proof}

\begin{Rem}
Note that the last part of the statement in Corollary \ref{action-cor} is a well known result which, for instance, follows from combining  \cite[Theorem 5.5]{sauer-betti-for-groupoids}  and \cite[Theorem 0.2]{cheeger-gromov}. 
\end{Rem}

Next we generalize the above result to the setting of amenable discrete measured groupoids;  $L^2$-Betti numbers for discrete measured groupoids were introduced by Sauer \cite{sauer-betti-for-groupoids} and generalize the notion for equivalence relations \cite{gaboriau-L2-for-eq-rel}. We recall only the basic definitions and notation here referring to \cite{sauer-betti-for-groupoids, renault-delaroche-amenable-groupoids-article} for the general theory of discrete measured groupoids. Recall that the object set $\ms G^0$ of a discrete measured groupoid $\ms G$ comes equipped with a  $\ms G$-invariant measure $\mu$ which gives rise to a measure $\mu_{\ms G}$ on $\ms G$ defined on a Borel subset $A\subseteq \ms G$ as
\beqn
\mu_{\ms G}(A):= \int_{\ms G^0} \card(s^{-1}(x)\cap A)\, d\mu(x).
\eeqn
Here $s\colon \ms G \to \ms G^0$ denotes the source map of the groupoid.
A measurable subset $K\subseteq \ms G$ is called \emph{bounded} if there exists an $M>0$ such that for almost every $x\in\ms G^0$ both $\card (s^{-1}(x) \cap K) < M$ and $\card (s^{-1}(x) \cap K^{-1}) < M$.  Let $\mb C\ms G$ be the groupoid ring of $\ms G$, defined as
\beqn
\mb C\ms G:=\{f\in L^\infty(\ms G,\mu_{\ms G})\,|\, \supp(f)\text{ is a bounded subset of }\ms G\}.
\eeqn
(this is well-defined independently of the choice of the representative of $f$). The convolution and involution on $\mb C\ms G$ are defined as
\beqn
(f\ast g)(\gamma) := \sum_{\substack{\gamma',\gamma''\in \ms G,\\\gamma'\gamma''=\gamma}} f(\gamma')g(\gamma''),
\eeqn
\beqn
f^*(\gamma) := \overline{f(\gamma^{-1})}.
\eeqn
This gives $\mb C\ms G$ the structure of a $\ast$-algebra which forms a strongly dense $\ast$-subalge\-bra in the groupoid von Neumann algebra $L\ms G$ \cite{renault-delaroche-amenable-groupoids-book}. For arbitrary subsets $K\subseteq \ms G$, $A\subseteq \ms G^0$, $B\subseteq \ms G^0$ we denote by $K_A^B$ the set $K\cap s^{-1}(A)\cap r^{-1}(B)$. Notice in particular that for every $x\in \ms G^{0}$ the set $\ms G^x_x$ is a group with respect to the composition. It is called the stabilizer group of $x$.
A discrete measured groupoid is called ergodic if for every set $A\subseteq \ms G^{0}$ of positive measure $r(s^{-1}(A))=\ms G^{0}$, and every discrete measured groupoid admits an ergodic decomposition into a direct integral of ergodic groupoids.
If $K\subseteq\ms G$ is a bounded measurable subset, then 
\[
 \mc F_K := \{f\in \mb C\ms G\,|\, \supp f \subseteq K\} 
\]
is a finitely generated projective right Hilbert $L^\infty(\ms G^{0})$-module when equipped with the  $L^\infty(\ms G^{0})$-valued inner product arising from the conditional expectation $E\colon L\ms{G}\to L^\infty(\ms{G}^0)$. This can be easily seen as follows:  take an $n\in\mb N$ such that almost every $x\in \ms G^0$ has at most $n$ preimages from $K$ with respect to $s$ and decompose $K$ as $K=K_1\sqcup K_2 \sqcup \dots\sqcup K_n$ in such a way that the map $s|_{K_i}$ is a measurable isomorphism onto its image $X_i\subseteq \ms G^0$ \cite[Lemma 3.1]{sauer-betti-for-groupoids}. Then $\mc F_K$ is isomorphic as a (pre\nobreakdash-)Hilbert \Cs-module over $L^\infty(\ms G^{0})$ to the direct sum $\oplus_{i=1}^n \mf{1}_{X_i} L^\infty(\ms G^{0})$ endowed with the standard inner product. This follows from the fact that for $f \in \mc F_{K_i}$ and $x\in X_i$ we have
\beqn
E(f^*f)(x) = (f^*\ast f)({\id}_x) = \sum_{\substack{\gamma',\gamma \in K_i,\\ \gamma'\circ \gamma = {\id}_x}} \overline{f\left({\gamma'}^{-1}\right)} f(\gamma)  = |f(\gamma(x))|^2,
\eeqn
where $\gamma(x)\in K_i$ is uniquely determined by the property $s(\gamma(x))=x$. Note in particular that $\dim_{L^\infty(\ms G^{0})}\mc F_K=\mu_{\ms G}(K)$. Moreover, a direct computation gives the following formula for $\varphi_{\mc F_K}$:
\beq\label{eq:rel-dim-for-groupoids}
\varphi_{\mc F_K}(T) = \frac{\sum_{i=1}^n \int_{X_i} E(T)\,d\mu}{\mu_{\ms G} (K)}= \frac{\sum_{i=1}^n \int_{X_i} E(T)\,d\mu}{\sum_{i=1}^n \mu(X_i)},\quad T\in L\ms G.
\eeq

The result now is as follows.

\begin{Cor}[Discrete measured amenable groupoids]\label{groupoid-example}
 Let $\ms{G}$ be an ame\-na\-ble dis\-crete measured groupoid for which the stabilizers $\ms G_x^x$ are finite for almost every $x\in \ms G^{0}$. Then the $L^2$-Betti numbers of $\ms G$ vanish in positive degrees:
\beqn
\beta_p^{(2)}(\ms G) := \dim_{L\ms{G}} \Tor_p^{\mb C\ms{G}}(L\ms{G},L^\infty(\ms G^0))=0. 
\eeqn
In particular, the $L^2$-Betti numbers of an amenable equivalence relation vanish.
\end{Cor}
Note that Corollary \ref{groupoid-example} is just a special case of \cite[Theorem 1.1]{sauer-thom}. In particular, the assumption that the stabilizers are finite is not important for the validity of the statement in Corollary \ref{groupoid-example}, but only needed in order for our methods to apply.

\begin{proof}
First of all, we observe that by \cite[Remark 1.7]{sauer-thom} it is sufficient to prove the statement for ergodic groupoids, and we will restrict ourselves to this case for the rest of the proof. Let $\ms R$ be the equivalence  relation  on $\ms G^{0}$ induced by $\ms G$; it can be considered as a quotient groupoid of $\ms G$ after dividing out all stabilizers and is known as the derived groupoid of $\ms G$. Since $\ms G$ is ergodic the same is true for $\ms R$. By amenability of $\ms G$, the quotient $\ms R$ is an amenable equivalence relation \cite[Theorem 2.11]{renault-delaroche-amenable-groupoids-article}, which is therefore hyperfinite \cite[Theorem XIII.4.10]{takesaki-3}. That is, there exists an increasing sequence of subrelations $\ms R_n$ each of which has finite orbits of cardinality $k_n$ and such that $\ms R\setminus \bigcup _n \ms R_n$ is a zero set. 
Let $\ms G_n$ be the subgroupoid of $\ms G$ generated by the stabilizers $\ms G_x^x$ and the lifts of the morphisms from $\ms R_n$ to $\ms G$. This subgroupoid is well-defined because two lifts of a morphism from $\ms R_n$ differ by an element in the stabilizer. Let $\mc A:= \varinjlim \mb C\ms G_n$ be the algebraic direct limit of the corresponding groupoid algebras. By construction of $\ms G_n$, $\ms G \setminus \bigcup_n \ms G_n$ is a measure zero subset of $\ms G$ and therefore $\mc A$ is an \emph{$L^\infty(\ms G^{0})$-rank dense} subalgebra of $\mb C\ms G$; i.e.~for each $f\in \mb C \ms G$ and for each $\eps >0$ there exists a projection $p\in L^\infty(\ms G^{(0)})$ with $\tau(p) > 1 - \eps$ such that $pf\in \mc A$. (for instance, one can take the characteristic function of a subset $X_f$ in $\ms G^{0}$ of measure bigger than $1-\eps$  with the property that $f$ is non-zero only on those morphisms from $\ms G \setminus \ms G_n$ whose target is in the complement of $X_f$). Furthermore, as $\mb C\ms G$ is dimension-compatible as an $L^{\infty}(\ms G^{0})$-bimodule \cite[Lemma 4.8]{sauer-betti-for-groupoids} the rank-density together with \cite[Theorem 4.11]{sauer-betti-for-groupoids} implies that
\beq\label{eq:l2-betti-numbers-hyperf-approx}
\beta_p^{(2)}(\ms G) := \dim_{L\ms G} \Tor^{\mb C\ms G}_p (L\ms G,L^\infty(\ms G^{0})) =\dim_{L\ms G} \Tor^{\mc A}_p (L\ms G,L^\infty(\ms G^{0})).
\eeq
Let us now consider $\mb C\ms G_n$ as a right $L^\infty(\ms G^{0})$-module. By construction, it is equal to $\mc F_{\ms G_n}$, where $\ms G_n$ is considered as a bounded subset of $\ms G$. Its dimension can be easily computed as 
$$\dim_{L^\infty(\ms G^{0})} \mb C \ms G_n = k_n\cdot \card \ms G_x^x,$$
where $\card \ms G_x^x$ is the cardinality of the stabilizer at $x \in \ms G^{0}$; it is an essentially constant function on $\ms G^{0}$ because of the ergodicity hypothesis. Moreover, \eqref{eq:rel-dim-for-groupoids} gives that
$$
\varphi_{\ms G_n}(T) = \tau(T),\quad T\in L\ms G.
$$
Now, putting $M:=L(\ms G)$ and considering $N:=L^\infty(\ms G^0)$ as a subalgebra of $\mc A$ via the inclusion $\ms G^0 \subseteq \ms G_n$, we infer that the subalgebra $\mc A$ has the strong F\o{}lner property. Therefore the inclusion $\mc A\subseteq M$ is dimension-flat and formula \eqref{eq:l2-betti-numbers-hyperf-approx} gives us the vanishing result for the  positive degree $L^2$-Betti numbers of $\ms G$.
\end{proof}

\begin{Ex}[The irrational rotation algebra]\label{irrational-rotation}
Let $\theta \in ]0,1[$ be an irrational number and recall  that the irrational rotation algebra $A_\theta$ is the universal $C^*$-algebra generated by two unitaries $u$ and $v$ subject to the relation $uv=e^{2\pi i\theta}vu$. 
The set 
\[
\mc A_\theta:=\left.\left\{ \sum_{p,q \in \ZZ } z_{p,q}u^{p}v^{q}\,\right|\, z_{p,q}\in \CC, z_{p,q}= 0 \text{ for all but finitely many } (p,q)\in \ZZ^2 \right\}
\]
forms a dense subalgebra in $A_\theta$ and carries  a natural filtering sequence of subspaces
\[
\mc F_n:= \left.\left\{ \sum_{-n\leqslant p +q\leqslant n} z_{p,q}u^{p}v^{q} \,\right|\, z_{p,q}\in \CC\right\}.
\]
Recall that $A_\theta$ has a unique trace which on $\mc A_\theta$ is given by $\tau\left(\sum_{p,q \in \ZZ } z_{p,q}u^{p}v^{q}\right)=z_{0,0}$.
We aim at proving that $\mc A_\theta $ has the strong F{\o}lner property relative to $N:=\CC$ and the enveloping von Neumann algebra $M_\theta:=A_{\theta}''\subseteq B(L^2(A_\theta,\tau))$. Let $\varps>0$ and $T_1,\dots, T_r\in \mc A_\theta$ be given. Then there exists an $m_0\in \NN$ such that $T_i(\mc F_n)\subseteq \mc F_{n+m_0}$ for each $i\in\{1,\dots, r\}$ and $n\in \NN$, and we now  define $\mc P_n := \mc F_{nm_0}$ and $\mc S_n:=\mc F_{(n-1)m_0}$. Using the fact that $H=L^2(A_\theta,\tau)$ has an orthonormal basis consisting of the unitaries
\[
\{ u^{p}v^{q}\mid p,q\in \ZZ \}
\] 
it is not difficult to see that $\varphi_{\mc P_n}(T)=\tau(T)$ for every $T\in M_{\theta}$ and furthermore it implies that $\dim_\CC{\mc P_n}=2(nm_0)^2+2nm_0+1$ from which it follows that $\frac{\dim_\CC \mc S_n}{\dim_\CC \mc P_n}\To 1$. Thus $(\mc P_n,\mc S_n)$ is a strong F{\o}lner sequence for $T_1,\dots,T_r$. From our results we therefore obtain that the Connes-Shlyakhtenko $L^2$-Betti numbers $\beta^{(2)}_p(\mc A_\theta,\tau)$ vanish for $p\geqslant 1$ as well as the well-known fact that the enveloping von Neumann algebra $M_\theta$ is hyperfinite.
\end{Ex}

\begin{Ex}
[UHF-algebras]\label{uhf-example}
Let  $A$ be  a UHF-algebra; i.e.~$A$ is a $C^*$-algebraic direct limit of a sequence $(\MM_{k(n)}(\CC),\alpha_n)$ of full matrix algebras. Denoting by $\mc A$ the algebraic direct limit we get a natural dense $\ast$-subalgebra in $A$ and  the tracial states on the matrix algebras $\MM_{k(n)}(\CC)$ give rise to a unique tracial state $\tau$ on $A$, for which the enveloping von Neumann algebra $M:=A''\subseteq B(L^2(A,\tau))$ is the hyperfinite II$_1$-factor. Denote by $\mc F_n$ the image of $\MM_{k(n)}(\CC)$ in $\mc A$. Then for given $T_1,\dots, T_r \in \mc A$ there exists an $n_0$ such that they are all in $\mc F_{n_0}$ and thus $T_i(\mc F_{n})\subseteq \mc F_{n}$ for all $i\in \{1,\dots, r\}$ and $n\geqslant n_0$. Putting $\mc P_n= \mc S_n:=\mc F_{n_0+n}$ we get a strong F{\o}lner sequence (relative to $\CC\subseteq \mc A\subseteq M)$ for $T_1,\dots, T_r$: the two first conditions in Proposition \ref{strong-foelner-prop} are obviously fulfilled and the a direct computation\footnote{This can be seen either using matrix units or a  collection of unitaries matrices forming an orthonormal basis for the Hilbert-Schmidt inner product.} 
shows that $\varphi_{\mc P_n}(T)=\tau(T)$. Thus the Connes-Shlyakhtenko $L^2$-Betti numbers of $(\mc A, \tau)$ vanish in positive degrees. 
\end{Ex}

\begin{Rem}
Note that by \cite[Corollary 3.2 \& Theorem 3.3]{piotr} neither the UHF-algebra nor the irrational rotation algebra is the algebra associated with a discrete group or quantum group; hence Example \ref{uhf-example} and Example \ref{irrational-rotation} provide honest new examples of dimension flat inclusions.
\end{Rem}

\begin{Rem} Generalizing Example \ref{uhf-example}, we may consider an arbitrary finite factor $N$ and an inductive system of the form $(\MM_{k(n)}(N),\alpha_n)$ and the inclusion of the algebraic direct limit $\mc A$ in the von Neumann algebraic direct limit $M$. Then, by exactly the same reasoning as in Example \ref{uhf-example}, $\mc A$ has the strong F{\o}lner property relative to the chain $N\subseteq \mc A\subseteq M$. Hence the inclusion $\mc A\subseteq M$ is dimension flat even though $\mc A$ can be ``far from amenable'' ($N$ might for instance be a II$_1$-factor with property (T)). Compare e.g.~with \cite[Conjecture 6.48]{Luck02} where it is conjectured that for $N=\CC$ and $\Gamma$ a discrete group, dimension flatness of the inclusion $\CC\Gamma\subseteq L\Gamma$ is equivalent to amenability of $\Gamma$.
\end{Rem}


\end{document}